\documentclass [12pt,english] {article}
\usepackage{amsthm}
\usepackage{amsmath}
\usepackage{amsfonts}
\usepackage{amssymb}
\usepackage[all]{xy}
\usepackage{color}   
\usepackage{tikz} \usetikzlibrary{arrows}

\usepackage{listings}
\usepackage{setspace}
\usepackage[font=small,labelfont=bf]{caption}
\singlespace
\def\TC{\protect\operatorname{TC}}

\def\zcl{\protect\operatorname{zcl}}

\def\hdim{\protect\operatorname{hdim}}

\newtheorem{ejem}{Example}
\newtheorem{defi}[ejem]{Definition}
\newtheorem{teo}[ejem]{Theorem}
\newtheorem{prop}[ejem]{Proposition}
\newtheorem{lema}[ejem]{Lemma}
\newtheorem{remark}[ejem]{Remark}
\newtheorem{coro}[ejem]{Corollary}

\numberwithin{ejem}{section}

\begin{document}
\title{Cup products and the higher topological complexity of configuration spaces of the circle with two anchored points}

\author{Teresa I. Hoekstra-Mendoza}

\date{\empty}

\maketitle
\begin{abstract}
    In this paper we show how to compute cup products in the anchored configuration space of the circle with two anchored points using discrete Morse theory. Knowing how to compute cup products allows us to obtain bounds for the (higher) topological complexity $TC_s$, which are sharp for a sufficiently large value of $s$.  
\end{abstract}

\section{Anchored configuration spaces}
Configuration spaces constitute a well-studied class of topological spaces. In usual configuration spaces, two or more particles are not allowed to occupy the same space, since we want to avoid collisions. In anchored configuration spaces collisions are allowed, but we  require that a certain pre-determined discrete set of points is always occupied by at least one particle.

Anchored configuration spaces are motivated by the problem of distributing $n$ unique resources to $m$ locations such that each location is never left without resource. In \cite{Ks} Kozlov studied this problem when the $m$ locations, to which the resources are distributed, are
connected by a tree. In this case, the anchored configuration space on $n$ points is called the Stirling complex $Str(T,n)$ where $T$ is a tree on $m$ vertices. Kozlov also proved that Stirling complexes are homotopy equivalent to a wedge of $n-m$-dimensional spheres and the amount of these spheres depends only on $m$ and $n$.

\

In \cite{K}, Kozlov gave an explicit combinatorial description of a basis for cohomology when there are two anchored points on a circle. We are going to use this cohomology basis to calculate cup products between them. Knowing how to compute cup products allows us to obatin bounds for the topological and higer topological complexity of a space.

\begin{defi}
Let $X$ be a non-empty topological space, let $S$ be a set of $m$ points in $X$, $m \geq  0$,
and let $n$ be an arbitrary positive integer. An anchored configuration space, denoted $\Sigma(X, S, n)$, is defined as the subspace of the direct product $X^n$, consisting of all tuples $(x_1 ,\dots, x_n )$, such that $S \subseteq \{x_1 , \dots , x_n \}.$
\end{defi}

\begin{remark} We have the following observations.
\begin{itemize}
\item When $m=0$, we simply have $\Sigma (X, \emptyset, n) = X^n$. Thus we may assume that $m>0$.
\item If $n<m$ we have $\Sigma(X, S, n) = \emptyset$.
\item If $n=m$, the space $\Sigma(X, S, n)$ is a collection of $n!$ points equipped with the discrete topology.
\end{itemize}
\end{remark}
This means we can assume that $n>m$. Denote by $K(X)$ the set of cells of a space $X.$

We want to focus on the case when $X$ is a circle, and the set $S$ has cardinality two.
In this case, $\Sigma(S^1,S,n)=\Omega_n$ has a cubical structure, where the cubes are tuples $(A,B,C,D)$ such that
\begin{itemize}
\item both $A$ and $C$ are non empty sets;
\item $A,B,C,D$ are pairwise disjoint;
\item $A\cup B \cup C \cup D =\{1,\dots, n\}$ and
\item the dimension of the cube is $|B|+|D|$.
\end{itemize}
The cubical boundary operator over $ \mathbb{Z}_2$ is given by $\partial (A,B,C,D) =$
$$ \sum_{x\in B} (A\cup x,B\setminus x,C,D) +(A,B\setminus x,C\cup x,D) + \sum_{x\in D} (A\cup x,B,C,D\setminus x)+(A,B,C\cup x,D\setminus x)  $$
We shall use the notation $\sigma = (A(\sigma),B(\sigma),C(\sigma),D(\sigma))$ as in Figure \ref{cubo}(left).
Throughout this paper we shall only consider coefficients in $\mathbb{Z}_2$.

\begin{figure}
\centering
\begin{tikzpicture}
\draw (2,0) circle [radius=2];
\node[circle, draw, scale=.4, fill=black] (a) at (0,0){};
\node (A) at (-.6,0){$A(\sigma)$};
\node[circle, draw, scale=.4, fill=black] (c) at (4,0){};
\node (C) at (4.6,0){$C(\sigma)$};
\node (b) at (2,2.5){$B(\sigma)$};
\node (d) at (2,-2.5){$D(\sigma)$};

\draw (9,0) circle [radius=2];
\node[circle, draw, scale=.4, fill=black] (a) at (7,0){};
\node (A) at (6.4,0){$\{1\}$};
\node[circle, draw, scale=.4, fill=black] (c) at (11,0){};
\node (C) at (11.6,0){$\{4\}$};
\node (b) at (9,2.5){$\{2,3\}$};
\node (d) at (9,-2.5){$\{5\}$};
\end{tikzpicture}
\caption{A cube $\sigma \in \Omega_n$, and an example of a cube in $\Omega_5$.}\label{cubo}
\end{figure}

\begin{ejem}
Take $n=5$  and consider the 3-cell $(1, \{2,3\}, 4, 5)$ as in Figure \ref{cubo} (right). The boundary of this cell is $(\{1,2\},3,4,5)+
(\{1,3\},2,4,5)+ (1,2,\{3,4\},5)+ (1,2,\{2,4\},5)+ (\{1,5\},\{2,3\},4, \emptyset) +( 1, \{2,3\}, \{4,5\}, \emptyset).$
\end{ejem}

However, this structure gives us a cubical complex structure but not a cubical set structure and we want to have a cubical set structure so we can apply the formula for cubical cup product given in \cite{KM}.

So consider the cycle graph on four vertices $C_4$, but assume that only two of the four vertices are anchored since we want to analyze the anchored configuration space of the circle with only two anchored points and assume these two are consecutive vertices.
We can construct a similar structure which is a cubical set structure as follows. Consider a spanning path in which the last two vertices are the anchored points. Then the cells are tuples $(x_1,y_1, \dots x_4,y_4)$ with the following properties:

Each $x_i$ and $y_i$ are pairwise disjoint sets for $1,\dots, 4$, such that $$x_3 \neq \emptyset \neq x_4\text{ and }\bigcup_{i=1}^4 (x_i \cup y_i)=\{1,\dots,n\}.$$ The dimension of the cell is $\sum_{i=1}^4 |y_i|.$

The cubical boundary operator over $\mathbb{Z}_2$ is given by $\partial(x_1,y_1, \dots, x_4,y_4)=$ $$\sum_{i=1}^4 \sum_{v\in x_i} ((\dots, x_i\setminus v, y_i\cup v,\dots)+(\dots, y_{i-1}\cup v, x_i\setminus v,\dots) $$
where we are considering $i$ modulo  $4.$
We shall denote this cubical set structure by $C\Omega_n$ and use the notation $\sigma=(x_1(\sigma),y_1(\sigma), \dots x_4(\sigma),y_4(\sigma)).$
\section{Discrete Morse theory}
A very power full tool for analyzing configuration spaces of graphs is discrete Morse theory, since it allows us to reduce the amount of cells in a complex while preserving the topological properties.
Assume throughout this section that $X$ is a regular complex.

\begin{defi}
Let $W$ be a collection of pairs of cells $(\sigma, \tau)$ such that $\sigma$ is a face of $\tau$ with $dim(\sigma)+1= dim(\tau).$ If each cell of $X$ appears as an entry of at most one pair of $W$ then $W$ is called a discrete vector field and
a cell $\sigma$ of $X$ is called:
\begin{itemize}
\item critical, provided it does not appear as an entry of any pair of $W$;
\item redundant, provided there is a cell $\sigma'$ such that $(\sigma,\sigma')\in W$;
\item collapsible, provided there is a cell $\sigma'$ such that $(\sigma',\sigma)\in W$.
\end{itemize}
  In other words, $W$ is a discrete vector field provided any cell of $X$ is of one and only one of the three types above.
\end{defi}
For a redundant cell $\tau$ of $X$, we shall denote by $W(\tau )$ the unique cell of $X$ with $(\tau, W(\tau))\in W$.
\begin{defi}
Let $W$ be a discrete vector field on $X$. A sequence of $k$-cells, $\tau_1,\ldots,\tau_n$ satisfying $\tau_i \neq \tau_{i+1}$ for $ i =1,\ldots,n-1$ is called an upper $W$-path of length $n$ if, for each $i=1,\ldots,n-1$, $\tau_i$ is redundant and $\tau_{i+1}$ is a face of $W(\tau_i)$.
Similarly, if $\tau_i$ is collapsible with $\tau_i=W(\sigma_i)$, and $\sigma_i$ is a face of $\tau_{i-1}$ for $i=2, \dots ,n$ then the sequence is called a lower $W$-path of length $n$.
 The $W$-path is closed if $\tau_1=\tau_n$. We say that $W$ is a gradient vector field provided it does not admit closed $W$-paths.
\end{defi}

We can also think of the gradient paths as directed paths in the Hasse diagram as follows.
Let $X$ denote a finite regular cell complex. The Hasse diagram of $X$, $H_X$ is a directed graph, where the vertices are the cells of $X$, and there exists an arrow from the vertex $v$ to the vertex $w$ if $w$ is a face of $v$, and $\text{dim} (v)= \text{dim} (w)+1$. Given a discrete vector field $W$ on $X$, the modified Hasse diagram $H_X(W)$ is the directed graph obtained from $H_X$ by reversing every arrow belonging to $W$. A $W$-path is a directed path in $H_X(W)$ which alternates reversed arrows with arrows in $H_X$.

We need to recall how gradient paths recover (co)homological information. In the rest of the section we assume $W$ is a gradient field on $X$.

Start by fixing an orientation on each cell of $X$ and, for cells $a^{(p)}\subset b^{(p+1)}$, consider the incidence number $\iota_{a,b}$ of $a$ and $b$, i.e. the coefficient ($\pm1$, since $X$ is regular) of $a$ in the expression of $\partial(b)$. Here $\partial$ is the boundary operator in the cellular chain complex $C_*(X)$. The Morse cochain complex $\mathcal{M}^*(X)$ is then defined to be the graded $R$-free\footnote{Cochain coefficients are taken in a ground ring $R$, as we are interested in cup products.} module generated in dimension $p\geq0$ by the duals\footnote{We omit the use of an asterisk for dual elements.} of the oriented critical cells $A^{(p)}$ of $X$. The definition of the Morse coboundary map in $\mathcal{M}^*(X)$ requires the concept of multiplicity of upper/lower paths. For a path of length two, multiplicity is given by
\begin{equation}\label{multiplicitydefinition}
\mu(a_0\nearrow b_1\searrow a_1)=-\iota_{a_0,b_1}\cdot\iota_{a_1,b_1}\mbox{ \ \ and \ \ }
\mu(c_0\searrow d_1\nearrow c_1)=-\iota_{d_1,c_0}\cdot\iota_{d_1,c_1},
\end{equation}
and, in the general case, it is defined to be a multiplicative function with respect to concatenation of  paths. The Morse coboundary is then defined by
\begin{equation}\label{morsecoboundary}
\partial(A^{(p)})=\sum_{B^{(p+1)}}\left(\sum_{b^{(p)}\subset B} \left(\iota_{b,B} \sum_{\gamma\in\overline{\Gamma}(b,A)}\mu(\gamma)\right)\right)\cdot B.
\end{equation}
In other words, the Morse theoretic incidence number of $A$ and $B$ is given by the number of  gradient paths $\overline{\gamma}$ from $B$ to $A$ counted with multiplicity $\mu(\overline{\gamma}):=\iota_{b,B}\cdot\mu(\gamma)$.

\section{Gradient paths}
In \cite{K} Kozlov gave a discrete gradient vector field $\mathcal{M}$ over $\Omega_n$ with trivial Morse boundary. Notice that $\{ \nu \in K(C\Omega_n):x_i(\nu)=y_i(\nu)=\emptyset\}= K(\Omega_n)$ as sets.
We are going to use Kozlov's discrete gradient vector field to create a gradient vector field $W$ on $C\Omega_n.$ 
Consider first the functions $\alpha_i, \beta : K(C\Omega_n) \rightarrow \{1,\dots, n\}$ for $ i=1,2,3$ given by $ \alpha_i (\sigma) =  \text{min} (x_i (\sigma) \cup y_i(\sigma))$ and $ \beta (\sigma)= \text{max} (y_3(\sigma) \cup x_4(\sigma))$.


\begin{defi} 
Define $W$ as follows:
\begin{itemize}
\item If $\alpha_i(\nu)\in x_i(\nu)$ for some $ i =1,2$ and if $i=2$ we have $x_1(\nu)=y_1(\nu)=\emptyset$ then $\nu $ is redundant and $W(\nu)$ is obtained from $\nu$ by moving $\alpha_i(\nu)$ from $x_i(\nu)$ to $y_i(\nu).$
\item If $\alpha_i(\nu)\in y_i(\nu)$ for some $ i =1,2$ and if $i=2$ we have $x_1(\nu)=y_1(\nu)=\emptyset$ then $\nu= W(\tau)$ is collapsible where $\tau$ is obtained from $\nu$ by moving $\alpha_i(\nu)$ from $y_i(\nu)$ to $x_i(\nu).$

   \item If $x_i(\sigma)=y_i(\sigma) = \emptyset$ for $i=1,2$ we take the gradient vector field $\mathcal{M}$ (taking $\alpha_3=\alpha$).
\end{itemize}
\end{defi}

It is easy to see that $W$ is well defined and is indeed a discrete vector field since a cell can either satisfy one or none of the above items. Moreover notice that the critical cells of $W$ are precisely the critical cells of $\mathcal{M}$ since a cell $\nu$, having either $y_i(\nu)$ or $x_i(\nu)$ non empty for some $i=1,2$, can never be critical.
So, recall which are the critical cells of $(\Omega_n,\mathcal{M})$

\begin{prop}{\cite{K}}
Let $$C_1= \{ \sigma \in \Omega_n : A(\sigma)= \alpha(\sigma), B(\sigma)= \emptyset, |C(\sigma)| \geq 2 \mbox{ and } \beta (\sigma) < \alpha (\sigma)\},$$
$$C_2= \{ \sigma \in \Omega_n : A(\sigma)= \alpha(\sigma), B(\sigma)= \emptyset,  \mbox{ and } C(\sigma)=\beta (\sigma) \}\text{ and }$$ 
$$C_3=\{ \sigma \in \Omega_n : A(\sigma)= \alpha(\sigma),B(\sigma) \neq \emptyset, C(\sigma)=\beta (\sigma) \mbox{ and } \alpha(\sigma)< B(\sigma)< \beta (\sigma) \} .$$
Then the set of critical cells of $\mathcal{M}$ is $C_1 \cup C_2 \cup C_3$.
\end{prop}

Notice that all critical cells $\sigma$ of dimension strictly smaller than $n-2$ have $y_3(\sigma) = \emptyset$. The following proposition gives us the amount of critical cubes in each dimension.

\begin{prop} \cite{K}
For $d \in \{0, 1, \dots, n-3\}$ the number of critical cubes of $W$ of dimension $d$ is $\binom{n}{d}$. The number of critical cubes of $W$ of dimension $n-2$ is $2^n+\binom{n-1}{2}-2.$
\end{prop}

\begin{prop}\label{nobi}
    Let $W(c)$ be a collapsible cell in $C\Omega_n$ such that $\alpha_i(c)\in y_i(c)$ for some $i=1,2$. Then there exists one unique non collapsible face of $W(c)$ which is obtained by moving $\alpha_i(c)$ from $y_i(c)$ to $x_{i+1}(c).$
\end{prop}

\begin{teo}
The discrete vector field $W$ is gradient.
\end{teo}
\begin{proof}
Assume that $W$ contains a cycle $C=c_1, \dots,c_k$ and notice that by definition of $W$, $y_4(c_1)=\dots =y_4(c_k).$
This implies that if $c_l=(\emptyset, \emptyset, \emptyset, \emptyset, x_3,y_3,x_4,y_4)$ for some $1\leq l \leq k$, thus for every cell $c_j$ with $j>l$, we have $x_i(c_j)=y_i(c_j)=\emptyset$ for $ i=1,2$. If every cell $c_j$ is such that $x_i(c_j)=y_i(c_j)=\emptyset$ we can obtain a cycle in $\mathcal{M}$ which is impossible.
By Proposition \ref{nobi}, there exists a unique way to obtain the redundant cell $c_{l+1}$ as a face of $W(c_l)$ for some $ 1 \leq l \leq k$.
 This means that if $v \in x_i(c_j)$ for some $ i=1,2$,$j\in \{1,\dots, k-1\}$  and $v \in x_{i+1}(c_{j+1})$ then $v \in x_l(c_m)$ for some $l>i$ and for every $j <m \leq k$. This means that $x_i(c_j)=x_i(c_l)$ for $i=1,2$ and for every $1 \leq j,l \leq k$ thus the cycle is trivial.
\end{proof}

\begin{lema}
    The Morse differentials in $(\Omega_n, \mathcal{M})$ and $(C\Omega_n, W)$ coincide.
\end{lema}
\begin{proof}
    Let $c=(\emptyset, \emptyset, \emptyset,\emptyset, x_3,y_3,x_4,y_4)$ be a critical cell. The only faces of $c$ which have one of its first four entries non empty, are the cells of the form $c=(v, \emptyset, \emptyset, \emptyset, x_3, y_3,x_4,y_4\setminus v)$.
    These cells are redundant and by Proposition \ref{nobi}, there exists a unique gradient path from $c$ to the cell  $e=(\emptyset, \emptyset, \emptyset,\emptyset, x_3\cup v,y_3,x_4,y_4)$. Notice that $e'=(x_3\cup v,y_3,x_4,y_4)$ is precisely the face of $(x_3,y_3x_4,y_4)$ in $\Omega_n.$ 
\end{proof}

\begin{teo} \cite{K} \label{teo0}
The Morse differential in $(\Omega_n, \mathcal{M})$ vanishes and therefore for each $m \geq 0$, a graded basis of $H^m(\Omega_n)$ is given by the cohomology classes of the duals of the critical $m$-cells.
\end{teo}
The following corollary allows us to obtain a basis for the cohomology of $C\Omega_n$ which we shall use to calculate cup products between basis elements.
\begin{coro}\label{0}
The Morse differential in $(C\Omega_n,W)$ vanishes and therefore for each $m \geq 0$, a graded basis of $H^m(\Omega_n)$ is given by the cohomology classes of the duals of the critical $m$-cells.
\end{coro}

\begin{coro}
For every $m\geq 0$, $H^m(\Omega_n)=H^m(C\Omega_n).$
\end{coro}

\
 
Since we are interested in critical cells, and gradient paths starting and ending at a critical cell, we define a 
 $\Lambda$-path to be an upper gradient path that ends at a critical cell.
Kozlov also proved the following lemma about $\Lambda$-paths which will be useful in the calculation of cup products.

\begin{lema}{\cite{K}} \label{lam}
\begin{itemize}
    \item Assume $\tau_1, \dots, \tau_q$ is a $\Lambda$-path. Then $B(\tau_i) =\emptyset$ for $1  \leq i \leq q$.
    \item Let $\sigma$ be a non-collapsible cube of $\Omega_n$. Then there exists a unique $\Lambda$-path starting at $\sigma$. Moreover, this path ends at the critical cell $\tau =(x, \emptyset, A(\sigma)\cup C(\sigma)\setminus\{x\}, D(\sigma))$ where $x= max\{A(\sigma) \cup C( \sigma)\}.$
\end{itemize}
\end{lema}
\begin{coro} \label{lam}
Let $c_1, \dots, c_q$ be a $\Lambda$-path in $(C\Omega_n, W).$ Then $y_i(c_j)=\emptyset$ for every $i=1,2,3,4 $ and $1\leq j \leq q.$ Moreover, if $\sigma$ is a non collapsible cube of $C\Omega_n$, then there exists a unique $\Lambda$-path starting at $\sigma$ which ends at the critical cell 
$$\tau =(\emptyset, \emptyset, \emptyset, \emptyset, x, \emptyset, x_3(\sigma)\cup x_4(\sigma)\setminus\{x\}, y_4(\sigma))$$ where $x= max\{x_3(\sigma) \cup x_4( \sigma)\}.$
\end{coro}

\section{Cup products}\label{cupprod}
We start this section by defining the cup product for cubical sets. To simplify notation, throughout this paper we shall omit the use of an asterisk for dual elements.

An elementary cube in $\mathbb{R}^k$ is a cartesian product $c=I_1\times\cdots\times I_k$ of intervals $I_i=[m_i,m_i+\epsilon_i]$, where $m_i\in\mathbb{Z}$ and $\epsilon_i\in\{0,1\}$. For simplicity, we write $[m]:=[m,m]$ for a degenerate interval. The standard product orientation of $c$ is determined by (a) the orientation (from smaller to larger endpoints) of the non degenerate intervals $I_{i_1},\ldots,I_{i_\ell}$ of $c$, and (b) the order $i_1<\cdots<i_\ell$, i.e. the order of factors in the cartesian product. Under these conditions, and for $1\leq r\leq \ell$, set
\begin{equation}\label{productboundaries}\begin{aligned}
\delta_{2r}(c)&=I_1\times\cdots\times I_{i_r-1}\times[m_{i_r}+1]\times I_{i_r+1}\times\cdots\times I_{k}, \\ \delta_{2r-1}(c)&=I_1\times\cdots\times I_{i_r-1}\times[m_{i_r}]\times I_{i_r+1}\times\cdots\times I_{k}.
\end{aligned}\end{equation}
Then, for a cubical set $X\subset\mathbb{R}^k$, i.e. a union of elementary cubes in $\mathbb{R}^k$, the boundary map in the oriented cubical chain complex $C_*(X)$ is determined by
\begin{equation}\label{productboundary}
\partial\left(c\right)=\sum_{r=1}^\ell(-1)^{r-1}\left(\rule{0mm}{4mm}\delta_{2r}(c)-\delta_{2r-1}(c)\right).
\end{equation}

\begin{figure}[h!]
$$
\begin{tikzpicture}[x=.4cm,y=.4cm]
\draw[->](0,0)--(3,0);\node[below] at (3,0) {\scriptsize${}+[0,1]\times [0]$};
\draw[->](3,0)--(6,0)--(6,3);\node[right] at (6,3) {\scriptsize${}+[1]\times [0,1]$};
\draw[->](6,3)--(6,6)--(3,6);\node[above] at (3,6) {\scriptsize${}-[0,1]\times [1]$};
\draw[->](3,6)--(0,6)--(0,3);\node[left] at (0,3) {\scriptsize${}-[0]\times [0,1]$};
\draw(0,3)--(0,0);
\node[left] at (0,0) {$(0,0)$};
\node[left] at (0,6) {$(0,1)$};
\node[right] at (6,6) {$(1,1)$};
\node[right] at (6,0) {$(1,0)$};
\end{tikzpicture}
$$
\label{orientedboundaryofsquare}
\caption{Oriented boundary of the square $[0,1]\times[0,1]$}
\end{figure}

Cup products in cubical cohomology are similar to their classic simplicial counterparts. At the oriented cubical cochain level, there is a cup product graded map $C^*(X)\times C^*(X)\to C^*(X)$ that is associative, $R$-bilinear and is described on basis elements as follows. Firstly, for intervals $[a,b]$ and $[a',b']$, let
$$[a,b]\cdot[a',b']:=\begin{cases} [a,b'], & \mbox{if $b=a'$ and either $a=b$ or $a'=b'$ (or both);} \\ 0, & \mbox{otherwise.}\end{cases}$$
Then, for elementary cubes $c=I_1\times\cdots\times I_k$ and $d=J_1\times\cdots\times J_k$ in $X$, the cubical cup product $c\cdot d$ of the corresponding basis elements\footnote{Recall that we shall omit the use of an asterisk for dual elements. The intended meaning will be clear from the context.} $c,d\in C^*(X)$ vanishes if either $I_i\cdot J_i=0$ for some $i\in\{1,\ldots,k\}$ or, if $(I_1\cdot J_1)\times\cdots\times(I_k\cdot J_k)$ is not a cube in $X$; otherwise $c\cdot d$ is the dual of the cube $(I_1\cdot J_1)\times\cdots\times(I_k\cdot J_k)$. 
Of course this definition involves a sign, but since in this paper we are considering coefficients in $\mathbb{Z}_2$ we may ignore them.

\

Translating this definition of cup products for cubical sets to our notation for $\Omega_n$ we obtain the following equation.
Let $(x_1,y_1,\dots, x_4,y_4)$ and $(z_1,w_1,\dots, z_4,w_4)$ be two cubes. Then their cup product as cubes is non zero if there exist sets $v_i$ for $i=1,2,3,4$ such that $x_i=v_i\cup w_i$ and $z_i=v_i\cup y_{i-1}$ for $ i=1,2,3,4$ putting $y_0=y_4$. Then 
$$(x_1,y_1,\dots, x_4,y_4)\cdot(z_1,w_1,\dots, z_4,w_4) = (v_1, y_1\cup w_1, \dots, v_4, y_4\cup w_4).  $$
Now we want to describe cup products in the Morse complex. So, assume $\sigma_1 $ and $\sigma_2$ are two critical cells. If there exist $\tau_1$ and $\tau_2$ such that there exist gradient paths $\tau_i \rightarrow \sigma_i$ for $i=1,2$ and $\tau_1 \cdot \tau_2 = \nu$ then the cup product $\sigma_1 \smile \sigma_2 = \nu_0$ where $\nu_0$ is a critical cell such that there exist a gradient path $\nu_0 \rightarrow \nu$ (see Figure \ref{copa}).

Here $\sigma_i$ and $\tau_i $ both have dimension $d_i$ for $i=1,2$ and $\nu $ and $\nu_0$ both have dimension $d_1 +d_2$.
\begin{figure}[h!]
\centering
\begin{tikzpicture}
\node (s1) at (0,2){$\sigma_1$};
\node (s2) at (0,0){$\sigma_2$};
\node (t1) at (2,2){$\tau_1$};
\node (t2) at (2,0){$\tau_2$};
\node (nu) at (3,1){$\tau_1 \cdot \tau_2=\nu$};
\node (n0) at (5,1){$\nu_0$};

\draw (t1)--(nu);
\draw (t2)--(nu);
\draw (n0)[->]--(nu);
\draw (t1)[->]--(s1);
\draw (t2)[->]--(s2);
\end{tikzpicture}
\caption{The cup product between critical cells.}\label{copa}
\end{figure}

\begin{lema}\label{inf}
Let $c$ be a collapsible cell and assume that $c=c_1\cdot, \dots,\cdot c_m$ where each $c_i$ is a cell of dimension $1$ for $1 \leq i \leq m.$ Then there exists a cell $c_j$ which is either collapsible, or there are no $\Lambda$-paths starting from $c_j.$
\end{lema}
\begin{proof}
Since $c$ is collapsible, $y_i(c)\neq \emptyset$ for some $i=1,2,3.$ Thus there exists $c_j$ such that $y_i(c_j)\neq \emptyset.$ If $i<3$ then $c_j$ is collapsible. If $i=3$, by Corollary \ref{lam}, there are no $\Lambda$-paths starting from $c_j.$
\end{proof}

\begin{lema}
The cup product of a critical cell with itself is always zero.
\end{lema}
\begin{proof}
Let $\sigma$ be a critical cell, clearly if $\text{dim}(\sigma) > \frac{n-2}{2}$ we have $\sigma \smile \sigma =0$.
So assume $\sigma = (\emptyset, \emptyset, \emptyset, \emptyset,x, \emptyset, Y,Z)$.  Consider $c_1=(x_1, \emptyset, x_2, \emptyset, x_3,\emptyset, x_4, y_4)$ and $c_2= (z_1, \emptyset, z_2, \emptyset, z_3,\emptyset, z_4, w_4)$ two cells such that the unique (recall Lemma \ref{lam}) $\Lambda$-path starting at them ends at $\sigma$. Then clearly $y_4 = w_4 = Z$. For the cells $c_1$ and $c_2$ to multiply each other we must have that $Z \subset x_1$ and $Z \subset z_3$ which is impossible since $c_1=(v_1,\emptyset,v_2,\emptyset,v_3, \emptyset, v_4 \cup Z, Z)$ and $c_2= (v_1\cup Z,\emptyset,v_2,\emptyset v_3, \emptyset, v_4, Z)$ for some sets $v_1, \dots, v_4$.  
\end{proof}

The following theorems tell us that we can multiply any two distinct critical 1-cells and their cup product is non zero. Moreover, this implies that any two distinct critical cells of dimensions $p$ and $q$ respectively have a non zero cup product as long as $p+q \leq n-2$.

\begin{teo}\label{prod}
Let $\sigma_i =(\emptyset,\emptyset,\emptyset,\emptyset, v_i, \emptyset, Z_i, w_i) $ be a critical 1-cell for $1 \leq i \leq j <n-2$. Then their cup product is the critical $j$-cell $(\emptyset,\emptyset,\emptyset,\emptyset ,a, \emptyset,C' , D)$ where $D=\{w_1, \dots, w_j\}$, $a= \text{max}\{\{1, \dots, n\}\setminus D\}\}$ and $C'= \{1, \dots, n\} \setminus (D\cup a)$.
\end{teo}
\begin{proof}
Since $v_i = \text{max}\{Z_i \cup v_i\}$, notice that we can reach the cell $\sigma_i$ with a gradient $\Lambda$-path starting at the cell $\tau_i = (B_i,\emptyset,E_i,\emptyset, A_i, \emptyset, C_i, w_i)$.

 We want to find conditions over $A_i,B_i,E_i,$ and $C_i$ such that $\prod_{i=1}^j \tau_i \neq 0$. This means that $B_i = B \cup \{w_1, \dots,w_{i-1} \}$, $ C_i= C\cup \{w_{i+1}, \dots, w_j\}$ for $1 \leq i  \leq j$, $E=E_i=E_k$ and $A=A_i=A_k$ for $1 \leq i,k \leq j$ for some sets $B,C$. 
 Then $\tau_1 \cdot \dots \cdot \tau_j=(B, \emptyset,E, \emptyset, A, \emptyset, C, \{w_1, \dots, w_j\})=\tau$. Since by Lemma \ref{inf} there are no lower gradient $\Lambda$-paths ending at $\tau,$ for the product $\sigma_1\smile \dots \smile \sigma_j$ to be non zero, $\tau$ must be critical. This means that $B=E=\emptyset$ and $A=max(A \cup C)$ thus $$\sigma_1 \smile \dots \smile \sigma_j = (\emptyset,\emptyset,\emptyset,\emptyset, a, \emptyset, C', D).$$
\end{proof}

\begin{coro}
Any critical cell of dimension $j <n-2$ is the product of $j$ critical 1-cells.
\end{coro}
\begin{proof}
Let $(\emptyset,\emptyset,\emptyset,\emptyset,A, \emptyset, C, D) $ be a critical $j$-cell, thus $A= \text{max}\{A \cup C\}$. Let $D=\{w_1,\dots, w_j\}$, then by Theorem \ref{prod}

$$(\emptyset,\emptyset,\emptyset,\emptyset,A,\emptyset,C,D)= \prod_{i=1}^j (v_i, \emptyset, (A\cup C\cup (D \setminus\{w_i\}))\setminus\{v_i\}, w_i)$$ where $v_i=max(A \cup C\cup (D\setminus\{w_i\}))$.
\end{proof}

\begin{teo}\label{n-2}
Let $\sigma_i =(\emptyset,\emptyset,\emptyset,\emptyset,v_i, \emptyset, Z_i, w_i) $ be a critical 1-cell for $1 \leq i \leq n-2$. Then  $\sigma_1 \smile \dots \smile \sigma_{n-2} =(\emptyset,\emptyset,\emptyset,\emptyset,A, \emptyset, C, D)+(\emptyset,\emptyset,\emptyset,\emptyset,C, \emptyset,A , D)$ where $D=\{w_1, \dots, w_{n-2}\}.$
\end{teo}

\begin{proof}
Since $v_i = \text{max}\{Z_i \cup v_i\}$, notice that we can reach the cell $\sigma_i$ with an upper gradient $\Lambda$-path starting at the cell $\tau_i = (B_i,\emptyset,E_i,\emptyset, A_i, \emptyset, C_i, w_i)$.

 We want to find conditions over $A_i,B_i,E_i,$ and $C_i$ such that $\prod_{i=1}^j \tau_i \neq 0$. This means that $B_i = B \cup \{w_1, \dots,w_{i-1} \}$, $ C_i= C\cup \{w_{i+1}, \dots, w_j\}$ for $1 \leq i  \leq n-2$, $E=E_i=E_k$ and $A=A_i=A_k$ for $1 \leq i,k \leq j$ for some sets $B,C$. 
 Then $\tau_1 \cdot \dots \cdot \tau_{n-2}=(B, \emptyset,E, \emptyset, A, \emptyset, C, \{w_1, \dots, w_{n-2}\})=\tau$. Since by Lemma \ref{inf} there are no lower gradient $\Lambda$-paths ending at $\tau,$ for the product $\sigma_1\smile \dots \smile \sigma_j$ to be non zero, $\tau$ must be critical. This means that $B=E=\emptyset$ and since $A$ and $C$ have no additional condition they need to satisfy, we can exchange their places to obtain a second critical cell as the product $\tau_1 \cdot \dots \cdot \tau_{n-2}$, thus $$\sigma_1 \smile \dots \smile \sigma_{n-2} = (\emptyset,\emptyset,\emptyset,\emptyset,A, \emptyset, C, D)+(\emptyset,\emptyset,\emptyset,\emptyset,C, \emptyset,A , D).$$
\end{proof}

\begin{coro}
Let $\sigma= (\emptyset,\emptyset,\emptyset,\emptyset,a, \emptyset, c, D)$ be a critical cell of dimension $n-2$. Then the sum $(\emptyset,\emptyset,\emptyset,\emptyset,a, \emptyset, c, D )+(\emptyset,\emptyset,\emptyset,\emptyset,c, \emptyset, a, D) $ can be factorized as the cup product of $n-2$ critical 1-cells. To be specific, we have $$ (\emptyset,\emptyset,\emptyset,\emptyset,a, \emptyset, c, D)+(\emptyset,\emptyset,\emptyset,\emptyset,c, \emptyset, a, D)= \prod_{i=1}^{n-2}(y_i, \emptyset, (\emptyset,\emptyset,\emptyset,\emptyset, a\cup c\cup (D \setminus\{w_i\}))\setminus\{v_i\}, w_i)$$ where $D=\{w_1,\dots, w_{n-2}\}$ and $v_i=max(a \cup c\cup (D\setminus\{w_i\}))$.
\end{coro}

\begin{remark}\label{obs}
\begin{enumerate}
\item Every critical cell of dimension smaller than $n-2$ has a unique factorization;
\item every critical cell $\sigma$ with $y_3(\sigma)=\emptyset$ and $\text{dim}(\sigma)=n-2$ appears in one unique product of $n-2$ critical cells of dimension one;
 \item a critical cell $\sigma$ with $y_3(\sigma)\neq\emptyset$ and $\text{dim}(\sigma)=n-2$ cannot be obtained as a product and does not appear as an element of any product.
\end{enumerate}
\end{remark}

\section{Topological complexity}
\begin{defi}
For any path-connected space $X$, the topological complexity of $X$, denoted by $TC(X)$, is the smallest integer $k \geq 0$ such that there is a cover of $X \times X$ by open sets $U_0, U_1, \dots, U_k$ and continuous sections  $s_i : U_i \rightarrow P(X)$ for $0 \leq i \leq k$. If there is no such $k$, set $TC(X)= \infty$.
\end{defi}

There is a cohomological bound for $TC(X)$. Let $R$ be a field and consider the cup product $$\smile: H^{\ast}(X;R) \otimes H^{\ast} (X;R)\rightarrow  H^{\ast} (X;R) .$$
Let $Z(X) \subset  H^{\ast} (X;R) \otimes  H^{\ast} (X;R)$ be the kernel of this homomorphism, called the ideal of zero divisors of  $H^{\ast} (X;R)$. The tensor product $ H^{\ast} (X;R)\otimes H^{\ast} (X;R)$ has a multiplication given by $(\alpha \otimes \beta)(\alpha'\otimes \beta') = (-1)^{|\alpha'|\cdot|\beta|} \alpha \alpha' \otimes \beta \beta'$, where $|x| = j$ if $x \in  H^{j} (X;R)$. The zero-divisors-cup-length of $ H^{\ast} (X;R)$ is the largest $i$ such that there exist elements $a_1, \dots, a_i \in Z(X)$ with $a_1 \cdot \dots, \cdot a_i \neq 0$.
Recall that we are taking coefficients over $\mathbb{Z}_2$.

\begin{teo}\cite{M}
Let $X$ be any CW complex, then $\text{zcl} (X) \leq TC(X) \leq 2 \text{dim} (X)$.
\end{teo}
We can use this last theorem, together with what we obtained in Section \ref{cupprod} to obtain the following bound.
\begin{teo}\label{cota}
For $n \geq 4$, we have the following bounds $n \leq TC(C\Omega_n) \leq 2n-4$.
\end{teo}
\begin{proof}
Since $C\Omega_n$ has dimension $n-2$ we know from Theorem \ref{cota} that $TC(C\Omega_n) \leq 2n-4$.
Let $\sigma_1, \dots, \sigma_n$ be the $n$ critical 1-cells. Let $\tau = \sigma_1 \smile \dots \smile \sigma_k$ and $ \nu = \sigma_{k+1} \smile \dots \smile \sigma_n$ for $k = \lfloor \frac{n}{2}\rfloor$.
Consider the zero divisors $\sigma_i \otimes 1 + 1 \otimes \sigma_i$ for $ 1 \leq i \leq n$.
We claim that the product
\begin{equation}\label{zd}
\prod_{i=1}^n (\sigma_i \otimes 1 + 1 \otimes \sigma_i)
\end{equation}
is non zero.
Product (\ref{zd})  contains the terms $\tau \otimes \nu$ and $ \nu \otimes \tau$. Since the cells $\nu$ and $\tau$ have a unique factorization, these two terms cannot appear elsewhere in the product (\ref{zd}) and thus it is non zero. By Theorem \ref{cota}, we have that $n\leq TC(C\Omega_n)-1$.
\end{proof}
\subsection{The $s$th topological complexity}
For $s\geq2$, the $s$th topological complexity of a path-connected space $X$, $\TC_s(X)$, is defined as the sectional category of the evaluation map $e_s\colon PX\to X^s$ which sends a (free) path on $X$, $\gamma\in PX$, to the $s$-tuple$$e_s(\gamma)=\left(\gamma\left(\frac{0}{s-1}\right),\gamma\left(\frac{1}{s-1}\right),\ldots, \gamma\left(\frac{s-1}{s-1}\right)\right).$$ 
 In other words, $\TC_s(X)+1$ stands for the smallest number of open sets covering $X^s$ in each of which $e_s$ admits a section. 

A standard estimate for the $s$th topological complexity of a space $X$, which we will use to obtain our bounds, is given by:
\begin{prop}\cite{B}\label{propcota}
For a $c$-connected space $X$ having the homotopy type of a CW complex, 
$$ \zcl_s(X)\leq\TC_s(X)\leq s\hdim(X)/(c+1).$$
\end{prop}

The notation $\hdim(X)$ stands for the (cellular) homotopy dimension of $X$, i.e.~the minimal dimension of a CW complex having the homotopy type of $X$.  
 The $s$th zero-divisor cup-length, $\zcl_s(X)$ is the largest non-negative integer $\ell$ such that there are classes $z_j\in H^*(X^s)$, each with trivial restriction under the iterated diagonal $\Delta_s \colon X \hookrightarrow X^s$, and such that the cup product $z_1\cdots z_\ell\in H^*(X^s)$ is non-zero. Each such class $z_i$ is called an $s$th zero-divisor for $X$. The ``zero-divisor'' terminology comes from the observation that the map induced in cohomology by $\Delta_s$ restricts to the $s$-fold tensor power $H^*(X)^{\otimes s}$ to yield the $s$-iterated cup product.

\begin{prop}\label{n}
Let $ s, n, j$ be positive integers such that $s \geq \lceil \frac{n}{j} \rceil$ for $ 2 \leq j \leq n$. Then there exist $s$ sets $V_1, \dots, V_s$ of cardinality $n-j$ such that $\bigcup_{i=1}^s V_i = \{1, \dots, n\} $ and $\bigcap_{i=1}^s V_i = \emptyset$.
\end{prop}
\begin{proof}
Assume for now that $s = \lceil \frac{n}{j} \rceil$ and let $U_1 = \{ 1 \dots, j\}$, $U_2=\{j+1, \dots, 2j\}$, $\dots,$ $U_s=\{ n-j+1, \dots, n\}$. Let $V_i = \{1, \dots, n\} \setminus U_i$ for $1 \leq i \leq s$. Then $|V_i|=n-j$,  $\bigcup_{i=1}^s V_i = \{1, \dots, n\} $ and $\bigcap_{i=1}^s V_i = \emptyset$.
If $s>\lceil \frac{n}{j} \rceil $ we can put $V_{\lceil \frac{n}{j} \rceil+j}=V_1$ for $1 \leq j \leq s - \lceil \frac{n}{j} \rceil$.
\end{proof}

\begin{teo}\label{bueno}
Let $s \geq \lceil \frac{n}{j} \rceil$ for $j \geq 2$. Then $s(n-j) \leq TC_s(C\Omega_n)$.
\end{teo}

\begin{proof}
Let $x_1, \dots, x_n$ be the $n$ critical cells of dimension one. By Proposition \ref{n} there exist $s$ sets $V_1,\dots, V_s$ of cardinality $n-j$ such that $\bigcup_{i=1}^s V_i = \{x_1, \dots, x_n\} $ and $\bigcap_{i=1}^s V_i = \emptyset$. Let $A_i$ be the cup product of the critical 1-cells in $V_i$ for $1 \leq i \leq s$. Then by theorems \ref{prod} and \ref{n-2}, $A_i$ is a critical $n-j$-cell if $j > 2$ or a sum of two critical $n-2$-cells if $j=2$.
We are going to define $s$th zero divisors as follows. Let $$\sigma_{i,l} = 1_1 \otimes 1_2 \otimes \dots \otimes 1_{i-1} \otimes x_l \otimes 1_{i+1}\otimes \dots \otimes 1_s + 1_1 \otimes 1_2 \otimes \dots \otimes 1_{i} \otimes x_l \otimes 1_{i+2}\otimes \dots \otimes 1_s \mbox{ if } x_l \in V_i$$  and $$\sigma_{s,l} = 1 \otimes 1 \otimes \dots \otimes 1 \otimes x_l + x_l \otimes 1 \otimes \dots \otimes 1 \mbox{ if } x_l \in V_s$$ for $1 \leq i \leq s-1$ and $1 \leq l \leq n$. We then have $s(n-j)$ zero divisors and we will show that their product is non zero.

Notice that the product 
\begin{equation}\label{sigma}
    \prod_{i=1}^s \prod_{x_l \in U_i} \sigma_{i,l}
\end{equation}
 contains the terms $A_1 \otimes A_2 \otimes \dots \otimes A_s + A_s\otimes A_1 \dots A_{s-1}$. By Remark \ref{obs} every $A_i$ for $1 \leq i \leq s$ has a unique factorization thus can not be obtained elsewhere in the product (\ref{sigma}) and hence do not cancel out.
This means that (\ref{sigma}) is non zero and by Proposition \ref{propcota} $s(n-j) \leq TC_s(C\Omega_n)$
\end{proof}

\begin{coro}
For $s \geq \lceil \frac{n}{2} \rceil$, the $s$th topological complexity $TC_s(C\Omega_n)$ is maximal.
\end{coro}
\begin{proof}
    By Theorem \ref{bueno} and Proposition \ref{propcota} we have $s(n-2) \leq TC_s(C\Omega_n) \leq s(n-2).$
\end{proof}

\section*{Acknowledgements}
The author would like to thank Jes\'us Gonz\'alez for noticing a mistake in the first version of this paper.

\end{document}